\documentclass{article}
\usepackage{fullpage}

\usepackage[utf8]{inputenc}

\usepackage{amssymb}
\usepackage{amsthm}
\usepackage{amsmath}
\usepackage{graphicx}

\newcommand{\eref}[1]{(\ref{#1})}

\newcommand{\imgext}[1]{#1.png} 

\begin{document}


\title{Hypercubic Self-Tilings}
\author{Ben Prather}


\begin{abstract}
In 1946 Fine and Niven posed problem E724, asking to demonstrate that every hypercube can be tiled by any number of hypercubic tiles larger than some value. This requires only basic number theory, but the problem of finding the smallest such number is much more involved.

For the square this is known to be 5, and the cube 47. No other values are known. This paper improves the bound for tesseracts from 808 to 733.
\end{abstract}

\maketitle
\section*{Introduction}
The story of this problem begins with American Mathematics Monthly problem E724.

\begin{quote}
Define a $D$-admissible number $t$ as one such that a $D$-dimensional hypercube may be subdivided into $t$ hypercubes. Prove that for each $D$ there exists an integer $A_D$ such that all integers exceeding $A_D$ are $D$-admissible\cite{FineNiven}.
\end{quote}

Several proofs were accepted and it was noted that finding the smallest such number, $h(D)$, is a much harder problem. In particular, the best reported bound on $h(3)$ was 54. Corner counting was used to rule out all cases except 54\cite{Scott}. 

In 1977 this was settled by the independent discoveries of cubic self-tiling with 54 cubes by Rychener and Zbinden\cite{Meier}. One was a flutist and the other an engineer, both from Switzerland where Hugo Hadwiger taught. About this time Hadwiger's name became associated with this problem despite no known published work on this problem. This is the reason we use $h(D)$.

Regardless, this established $h(2)=5$ and $h(3)=47$. Exact values are not known for any other dimension. 

In 1991 this problem appeared in the book Unsolved Problems in Geometry, volume II of Springer's Unsovled Problems in Intuitive Mathematics. In 1998 Hudelson improved the bound reported there from $h(4)\leq 853$ to $h(4)\leq808$\cite{Hudelson}. Hudelson's work was included in the 2001 edition of the above titled book.

In 2003 Erich Friedman posted a power point slide to his website featuring several open problem, attributing this one to Hugo Hadwiger and repeating the bound of $853$. It was here that I first found this problem in 2010. 

Within a couple weeks I had improved the bound to 838, unaware of Hudelson's work. Within a few months this was reduced to 763. In 2014 I tried, unsuccessfully, to publish this result. This lead to several new tilings, and the current bound of $h(4)\leq 733$.

\section{Solution to E724}
Observe that if the self tiling has more than 1 tile there must be a tile in every corner. Thus there must be at least $2^D$ tiles. Also, for every $D$ there is a tiling of $2^D$ tiles; simply cut in half in each direction.

Nesting this tiling into any other tiling replaces one tile with $2^D$ tiles, for a net gain of $2^D-1$ tiles. This reduces the problem to finding the smallest tilings mod $2^D-1$. 

Actually, for E724 all that is needed is any tiling with $t$ tiles such that $GCD(t-1,2^D-1)=1$. The tiling cutting each direction into $2^D-1$ equal portions provides $t=(2^D-1)^D$, which clearly works. This produces a bound of $h(D)\leq 2^{D^2+D}$.

\section{Improving the bound}
In 1976 Erdős improved this bound to $h(D)\leq e(2D)^D$. The best known asymptotic bound of $e(2D)^{D-1}$ is due to Hudelson\cite{Hudelson}. In most cases we can do much better. In particular, Hudelson showed that if $\gcd(2^D-1,k^D-1)=1$ then the bound is $(2k)^{D-1}$\cite{Hudelson}. 

Fermat's little theorem states that if $D+1$ is prime, then $n^{D+1}\equiv n\bmod D+1$. Thus $n^D\equiv1$, unless $n$ is a multiple of $D+1$, where $n^D\equiv0$.
This means that if $D+1$ is prime, then $D+1$ divides $k^D-1$ for all $k<D+1$ and the bound becomes $(2(D+1))^{D-1}<e(2D)^{D-1}$.

Since there are infinitely many primes this bound will occur infinitely many times and this is the best Hudleson's methods achieve.

Experience suggests that reducing the problem to cases mod $D+1$ happens rather quickly. Taking a $2^D$ tiling and replacing up to $D$ of the tiles with $(D+1)^D$ subtiles each gives an example of every case. In general, there should be plenty of room to merge tiles to produce each needed case plus a few copies of the original $2^D$ tiling.

This suggests a bound of $(D+1)^{D+1}<eD^{D+1}$, a possible asymptotic improvement by a factor of $2^D/D$. It should be noted that Hudelson's bounds are tighter for $D<7$, or if $k<D/2$. We now proceed to examine individual bounds for $D\leq7$. 

\section{$D=2$}
Hudelson's bound is $h(D)<11$. We need to find the best tilings mod $2^2-1=3$. $t=1$ is trivial. It is useful to have a notation to talk about a hypercube of edge $n$ in dimension $D$, say $H^D_n$. 

\begin{figure}[ht]
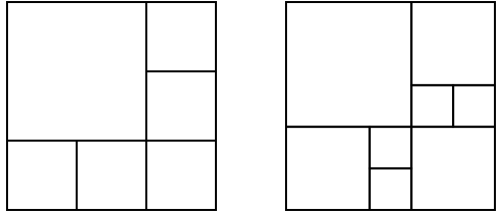

\centering
\includegraphics[width=3.0cm,height=3.0cm]{\imgext{C3}}\qquad
\includegraphics[width=3.0cm,height=3.0cm]{\imgext{E_5_3}}\qquad
\caption{Images of tilings in $D=2$ with $t=6$ and $t=8$ respectively.}
\label{fig:2DBest}
\end{figure}

$H^2_3=H^2_2+5\cdot H^2_1$ provides a tiling with $6\equiv 0$ tiles. Finally, $H^2_5=H^2_3+3\cdot H^2_2+4\cdot H^2_1$ provides a tiling with $8\equiv 2$ tiles. Note that $H^2_3\neq 2\cdot H^2_2+H^2_1$, despite $9=2\cdot4+1$. Thus volume alone is not sufficient to guarantee a tiling. The required tilings are shown in Figure \ref{fig:2DBest}. 

This set of tilings generates tilings for all $t$ other than 2,3 and 5. Note for example that $t=7$ can be attained by nesting two $t=4$ tilings. $t=2$ and $t=3$ are impossible since each corner must be in a tile, and if two corners share a tile it must be the trivial tiling.

Consider the largest corner square. If it is smaller than half the edge length there must be gaps on all edges. If it is larger than half the edge lengths there must be gaps on each edge adjacent to the corner opposite the larger corner tile. If it is exactly half, either we have $t=4$ or some corner has a tile less than half. In the latter case the edges of the corner with the smaller tile must have two gaps.

But this is all cases, so $t=5$ is not possible. This establishes that $h(2)=5$.

\section{$D=3$}
Hudelson's bound is $h(3)<98$. We need to find the best tiling mod $2^3-1=7$. 1 is trivial. $H^3_3=H^3_2+19\cdot H^3_1$ provides $t=20\equiv 6$. Nesting this tiling inside itself yields $t=39\equiv 4$. Similarly $H^3_4=H^3_3+37\cdot H^3_1$ provides $t=38\equiv 3$.

The cases known prior to 1977 and the lone improvement are:
\begin{itemize}
\item $H^6_3=4\cdot H^3_3+9\cdot H^3_2+36\cdot H^3_1$ yields $t=49\equiv 0$.
\item $H^3_6=5\cdot H^3_3+5\cdot H^3_2+41\cdot H^3_1$ yields $t=51\equiv 2$.
\item $H^3_6=3\cdot H^3_3+11\cdot H^3_2+47\cdot H^3_1$ yields $t=61\equiv 5$.
\item $H^3_8=6\cdot H^3_4+2\cdot H^3_3+8\cdot H^3_2+42\cdot H^3_1$ yields $t=54\equiv 5$.
\end{itemize}

\begin{figure}[ht]
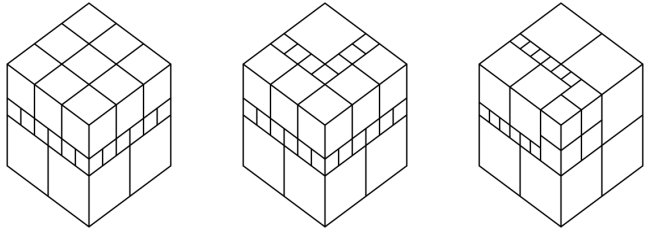

\centering
\includegraphics[width=2.43cm,height=3.24cm]{\imgext{M2_3_2_4}}\qquad
\includegraphics[width=2.43cm,height=3.24cm]{\imgext{M2_3_2_3}}\qquad
\includegraphics[width=2.43cm,height=3.24cm]{\imgext{M2_4_3_2}}
\caption{Tilings for $D=3$ for $t=49$, $t=51$, and $t=54$ respectively. The hidden corner contains another of the larger tiles in each. Several unit tiles are hidden internally.}
\label{fig:3DBest}
\end{figure}

Altogether this establishes $h(3)\leq 47$. The tilings for $t\in\{49,51,54\}$ are shown in Figure \ref{fig:3DBest}.

The corner and edge counting methods demonstrated above for $D=2$ were used to rule out all remaining values for $D=3$. Thus $h(3)=47$.

It is now useful to introduce a notation for some common families of tilings. Let $T_n$ represents an $H^D_n$ tiled with an $H^D_{n-1}$ in one corner and $H^D_1$ elsewhere. $T_{6-d}$ and $T_{8-d}$ are used to indicate a 2 tiling where $d$ corners are tiled further following the patterns shown above, producing an $H^D_6$ or $H^D_8$ respectively. For example, the tilings used for $D=3$ are $T_1$, $T_3$, $T_4$, $T_{6-2}$, $T_{6-4}$, $T_{6-3}$, $T_{8-2}$. Note that $T_{6-2}$ has the same number of tiles as 3 applied twice for any $D$.

The corner counting methods lead to considering exactly the geometries implemented by the $T_{6-i}$ and $T_{8-i}$ tilings. This suggests that enumerating these tilings would allow us to establish $h(D)$ if $GCD(2^D-1,3^D-1)=1$. It is not known whether this happens for infinitely many $D$. It does happen for $D=5$ and $D=7$.

\section{D=5 or 7}
Hudelson's bound is $h(5)<27183$ (3523 if we use $k=3$). We need to find the best tilings mod $2^5-1=31$. 

The lowest values found using these tilings for $D=5$ are given in Table \ref{tab:5D}.
\begin{table}[ht]
\caption{5D Tiling Types}
\label{tab:5D}
\centering
\begin{tabular}{ccc|ccc|ccc}
$t\mod31$& Tiling     & $t-31$ &$t\mod31$& Tiling & $t-31$  &$t\mod31$& Tiling & $t-31$ \\
\hline
0  & $T_{8-5}$  & 1705   & 11 & $T_{6-19}$ & 1809   & 22  & $T_{6-12}$  & 1541 \\
1  & $T_1$      & 1      & 12 & $T_{6-24}$ & 1655   & 23  & $T_{6-17}$  & 1573 \\
2  & $T_{6-5}$  & 901    & 13 & $T_{8-2}$  & 1098   & 24  & $T_{6-22}$  & 1915 \\
3  & $T_{6-10}$ & 1429   & 14  & $T_{6-15}$ & 1647  & 25  & $T_{8-4}$   & 1420 \\
4  & $T_{6-15}$ & 1647   & 15 & $T_{6-8}$  & 1131   & 26  & $T_3$     & 212  \\
5  & $T_{6-20}$ & 1803   & 16  & $T_{6-13}$  & 1659 & 27  & $T_{6-6}$   & 1019\\
6  & $T_{8-6}$  & 1742   & 17  & $T_{6-18}$  & 1691 & 28  & $T_{6-11}$  & 1547\\
7  & $T_4$      & 782    & 18  & $T_{8-8}$   & 1692 & 29  & $T_{6-16}$  & 1393\\
8  & $T_{6-4}$  & 721    & 19  & $T_{8-3}$   & 1383 & 30  & $T_{6-21}$  & 1921\\
9  & $T_{6-9}$  & 1311   & 20  & $T_{6-2}$   & 423  &&&\\
10 & $T_{6-14}$ & 1653   & 21  & $T_{6-7}$   & 1137 &&&\\

\hline
\end{tabular}
\end{table}
The worst case it $t=1921\equiv 30$. Subtracting 31 yields $h(5)\leq 1890$. Hudelson established this bound by tabulating the $T_{6-i}$ tilings alone\cite{Hudelson}. Several individual cases are improved from the Hudelson tabulation. 

Similarly $h(7)\leq67374$. It is likely that corner and edge counting can be used to rule out all remaining cases for both of these values.

\section{Intermission}
Before moving on to $D=4$, lets count some arrangements in $D=2$ and $D=3$ that will aide our counting.

Consider a tiling of a $H^D_5$ similar to the tiling of a square with 8 squares provided above. We thus put an $H^D_3$ in one corner, an $H^D_2$ in the remaining corners and fill the spaces with $H^D_1$. The following are of particular interest.

\begin{align}
D=2\qquad&H^2_5=H^2_3+3\cdot H^2_2+4\cdot H^2_1 \label{eq:H5_2}\qquad t=8\\
D=3\qquad&H^3_5=H^3_3+7\cdot H^3_2+42\cdot H^3_1 \label{eq:H5_3}\qquad t=50\\
D=4\qquad&H^4_5=H^4_3+15\cdot H^4_2+304\cdot H^4_1 \label{eq:H5_4}\qquad t=320
\end{align}

\begin{figure}[ht]
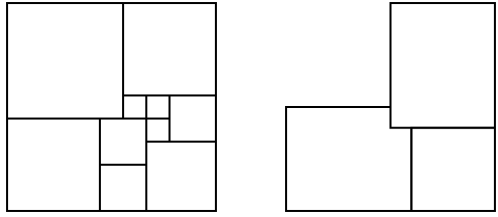

\centering
\includegraphics[width=3.0cm,height=3.0cm]{\imgext{S3_9_5_2}}\qquad
\includegraphics[width=3.0cm,height=3.0cm]{\imgext{t_698}}
\caption{A tiling of $H^2_9$ and $H^2_{10}-H^2_5$.}
\label{fig:H2_9_3}
\end{figure}

Figure \ref{fig:H2_9_3} demonstrates:
\begin{align}
H^2_9&=H^2_5+2\cdot H^2_4+H^2_3+3\cdot H^2_2+3\cdot H^2_1 \label{eq:H2_9}\\
H^2_{10}-H^2_5&=H^2_4+6\times5 + (6\times5-H^2_1)\label{eq:H2_10}
\end{align}

\begin{figure}[ht]
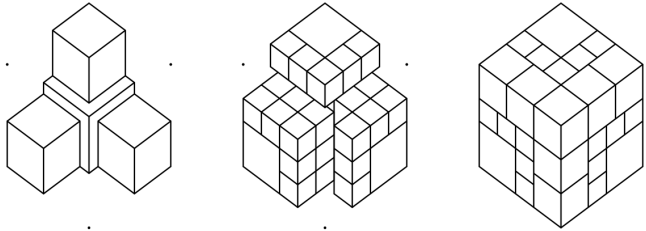

\centering
\includegraphics[width=2.43cm,height=3.24cm]{\imgext{S3_9_3a}}\qquad
\includegraphics[width=2.43cm,height=3.24cm]{\imgext{S3_9_3b}}\qquad
\includegraphics[width=2.43cm,height=3.24cm]{\imgext{S3_9_3c}}
\caption{Images showing layers of the tiling of an $H_9^3$. From left to right, the first pane shows the placement of the $H_5^3$ and three $H_4^3$. The second pane shows the placement of the 22 $H_2^3$. The final pane shows the placement of the 7 $H_3$. The 47 $H_1^3$ are hidden internally.}
\label{fig:H3_9_3}
\end{figure}

Figure \ref{fig:H3_9_3} demonstrates:
\begin{equation}\label{eq:H3_9}
H^3_9=H^3_5+3\cdot H^3_4+7\cdot H^2_3+22\cdot H^3_2+47\cdot H^3_1
\end{equation}

We now return to the regularly scheduled topic.

\section{$D=4$}
Hudelson's theoretical bound is $h(4)<1392$. We need to find the best tilings mod $2^4-1=15$. 1 is trivial. $T_3$ has $t=66$. Twice we get $t=131$. \eref{eq:H5_4} gives a tiling with $t=320$ and $T_5$ yields $t=370$. Nesting \eref{eq:H5_4} and $T_3$ gives $t=385$. These are the best known cases for $t\equiv 1$ or $t\equiv 0$ mod 5.

We now consider the remaining 9 cases separately, but grouped into equivalence classes mod 5.

\subsection{$t\equiv 2$}
These cases were not explored until my first attempt at publication since $T_6$ nested with $T_3$ is sufficient to establish the bound.
\subsubsection{$t=552$}
An $H_9^4$ can be tiled by 57 $H_3^4$ to leave a hyperbox of size $9\times6^3$. Partition this hyperbox into $5\times6^3$ and $4\times6^3$ regions. The first region can be tiled by an $H_5^4$ and 455 $H_1^4$, while the second is tiled by an $H_4^4$ and 38 $H_2^4$.

Altogether this totals 1 $H_5^4$, 1 $H_4^4$, 57 $H_3^4$, 38 $H_2^4$ and 455 $H_1^4$, producing $t=552$.

\subsubsection{$t=502$}
An $H_9^4$ can be tiled by 45 $H_3^4$ to leave a hyperbox of size $9^2\times6^2$. 

Using \eref{eq:H2_9}, the $9^2\times6^2$ hyperbox is now separated into a $5^2\times6^2$, two $4^2\times6^2$, a $3^2\times6^2$, three $2^2\times6^2$ and three $1^2\times6^2$ hyperboxes.

The $3^2\times6^2$, $2^2\times6^2$ and $1^2\times6^2$ hyperboxes are tiled by 4 $H_3^4$, 9 $H_2^4$ and 36 $H_1^4$ respectively. The $5^2\times6^2$ is tiled by an $H_5^4$ and 275 $H_1^4$. The two $4^2\times6^2$ regions are tiled by an $H_4^4$ and 20 $H_2^4$.

Altogether, this is 1 $H_5^4$, 2 $H_4^4$, 49 $H_3^4$, 67 $H_2^4$ and 383 $H_1^4$, producing $t=502$.

\subsubsection{$t=542$}
An $H_9^4$ can be tiled by 27 $H_3^4$ to leave a $9^3\times6$ hyperbox. This can be tiled by \eref{eq:H3_9} stacked 6 high. Each $H_3^3$, $H_2^3$ and $H_1^3$ in \eref{eq:H3_9} represents 14 $H_3^4$, 66 $H_2^4$ and 282 $H_1^4$ respectively. The three $4^3\times6$ hyperboxes can be tiled by an $H_4^4$ and 8 $H_2^8$ each. Likewise, the $5^3\times6$ hyperbox can be tiled by an $H_5^4$ and 125 $H_1^4$.

The grand total is now an $H_5^4$, 3 $H_4^4$, 41 $H_3^4$, 90 $H_2^4$ and 407 $H_1^4$, producing $t=542$. 

\subsection{$t\equiv 4$}
\subsubsection{$t=499$}
This tiling of an $H_{10}^4$ begins with 14 $H_5^4$, leaving a $10\times5^3$ hyperbox to be tiled. This can be partitioned using \eref{eq:H5_3} stacked as full as possible. A layer of $H_1^4$ will be required to finish the $H_3^4$ stack.

Altogether this requires 14 $H_5^4$, 3 $H_3^4$, 35 $H_2^4$,  and 447 $H_1^4$. This produces $t=499$.

\subsubsection{$t=534$}
This tiling is similar to the 499 tiling, except an $H_4^4$ replaces an $H_3^4$. The $10\times5^3$ region can then be separated into $6\times5^3$ and $4\times5^3$ regions.

The $4\times5^3$ hyperbox can be tiled by an $H_4^4$ and 244 $H_1^4$.
The $6\times5^3$ can be tiled by stacking \eref{eq:H5_3}. This requires 2 $H_3^4$, $21$ $H_2^4$ and 252 $H_1^4$.

Altogether this is a tiling with 14 $H_5^4$, 1 $H_4^4$, 2 $H_3^4$, 21 $H_2^4$ and 496 $H_1^4$. This produces $t=534$. 

\subsubsection{$t=584$}
This tiling is similar to the 534 tiling, except it starts with a $T_3$ instead of $T_2$. This adds $T_3-T_2=50$ tiles to the above case, producing $t=584$. 

This trick was found during the initial write up, improving my prior best bound. It was almost neglected in this case, until after several readings of the nearly finished paper. That makes this the last case found.

\subsection{$t\equiv 3$}
\subsubsection{$t=693$}
This tiling begins with an $H_{10}^4$, and tiles it with 13 $H_5^4$, leaving an L shaped $(10^2-5^2)\times5^2$ region. The $5^2$ component is then separated using \eref{eq:H5_2}.

The $(10^2-5^2)\times3^2$ region will fit 5 $H_3^4$ along the outer edges, plus 5 $H_2^4$ inside them. The rest of the space is tiled with 190 $H_1^4$. The three $(10^2-5^2)\times2^2$ regions can be tiled with 16 $H_2^4$ each, as shown by nesting a $T_6$ inside a $T_4$ inside a $T_5$ for $D=2$. This produces 16 $H_2^4$, with 144 $H_1^4$ left over in each region.

Altogether, this totals 13 $H_5^4$, 5 $H_3^4$, 53 $H_2^4$, and 622 $H_1^4$, producing $t=693$.

\subsubsection{$t=698$}
This tiling begins with leaving an l shaped region in an $H_{10}^4$ as above. We can separate the $(10^2-5^2)$ region using \eref{eq:H2_10}. 

The $4^2\times5^2$ region can be tiled by 1 $H_4^4$ and 144 $H_1^4$. The $6\times5^3$ region can be tiled by stacks of \eref{eq:H5_3}, producing 2 $H_3^4$, 21 $H_2^4$ and 252 $H_1^4$.

The $(6\times5-1^2)\times5^2$ region is similar except that it is interrupted by the blockage. This region can be partitioned using \eref{eq:H5_2} to produce a $(6\times5-1^2)\times3^2$, three $(6\times5-1^2)\times2^2$ regions and four $(6\times5-1^2)\times1^2$ regions. 

The $(6\times5-1^2)\times3^2$ region can be tiled by 2 $H_3^4$, 2 $H_2^4$ and 81 $H_1^4$. The $(6\times5-1^2)\times2^2$ regions can be tiled by 6 $H_2^4$ and 20 $H_1^4$ each. The $(6\times5-1^2)\times1^2$ regions can be tiled by 29 $H_1^4$ each. 

Altogether this yields 13 $H_5^4$, 1 $H_4^4$, 4 $H_3^4$, 41 $H_2^4$, and 639 $H_1^4$, producing $t=698$. 

\subsubsection{$t=748$}
Applying the $T_3$ trick used for $t=584$ to the $t=698$ tiling produces $t=748$. Since this is the worst case, we have $h(4)\leq 733$, the next lower number mod 15.

\subsection{Summary}
Hudelson established his bound of $h(4)\leq808$ using the tilings $T_2$, $T_3$, $T_5$, $T_6$, the tiling above for $t=678$ and a tiling with $t=619$\cite{Hudelson}. The improvements from these numbers are given in the Table \ref{tab:4D}.

\begin{table}[ht]
\caption{4D Tiling Types}
\label{tab:4D}
\centering
\begin{tabular}{cccc}
$t\mod15$ & New & Old & $\Delta$\\
\hline
0  & 435 & 435 & 0   \\
1  & 1   &  1  & 0   \\  
2  & 542 & 737 & 195 \\
3  & 693 & 693 & 0   \\
4  & 499 & 634 & 135 \\ 
5  & 320 & 500 & 180 \\
6  & 66  &  66 & 0   \\
7  & 502 & 802 & 300 \\
8  & 698 & 758 & 60  \\
9  & 534 & 699 & 165 \\
10 & 370 & 370 & 0   \\
11 & 131 & 131 & 0   \\
12 & 552 & 672 & 120 \\
13 & 748 & 823 & 75 \\
14 & 584 & 764 & 180 \\
\hline
\end{tabular}
\end{table}

Improvements are made for 9 cases, 6 by over 100. The largest improvement is 300 for the case $n\equiv7\mod 15$. Both the old and new bound are located at $n\equiv13\mod 15$, with the bound improved by 75 to $h(4)\leq733$. 

\section{D=6}
Hudelson's theoretical bound is $h(6)<676396$. We need to find the best tilings mod $2^6-1=63$. In a first pass $T_2$, $T_4$ and $T_{6-i}$ can be used to reduce the problem to finding the best cases mod 7.

Partition an $H^6_{14}$ like a $T_2$, sans 2 or 4 $H^6_7$ and putting 2 or 4 $H^6_5$ packed into one of the open corners. Filling the rest with $H^6_4$ , $H^6_2$  or $H^6_1$, space permitting to create $T_{14-2}$ and $T_{14-4}$. Combined with $T_8$, these tilings establishes $h(6)\leq246963$. 

This bound is very loose. In particular, no $H_3^6$ were used to fill the $H^6_{14}$. Introducing these strategically should lower this bound. Significant improvements to specific cases are also likely from other irregular tilings, in particular for the cases covered using $T_8$.

\section{Discussion}
For $D=4$ improvements have been found for several cases and the bound from $h(4)\leq808$ to $h(4)\leq733$. It is natural to conjecture that no further improvements can be made, setting $h(4)=733$. Alternatively, it is possible that efficient geometries have been missed. Tilings of an $H_{12}^4$ are strong candidates due to the large number of factors of 12. Fermat's little theorem requires these must have 2, 7 or 12 $H_5^4$.

For $D$ where $\gcd(2^D-1,3^D-1)=1$ the tilings $T_2$, $T_3$, $T_4$, $T_{6-i}$ and $T_{8-i}$ produce the best known tilings for each case. It is reasonable to conjecture that this must be the case. This would set $h(5)=1890$ and $h(7)=67374$.

Further, $T_4$ and $T_{8-i}$ appear to not improve the worst case, except for $D=3$. Thus it is reasonable to further conjecture that $h(D)=b$, where $b$ is the bound determined using $T_2$, $T_3$ and $T_{6-i}$ alone for $D>3$; as explored by Hudelson up to $D<25$\cite{Hudelson}.

The case of $D=6$ is entirely different, as even the tilings outlined here are left unexplored. Improvements to the bound of $h(6)\leq246963$ appear possible in a straight forward manner.

$h(D)$ is not expected to be monotonic, though a lower bounds for $h(6)$ would be needed to prove $h(7)\leq (6)$. The general problem of establishing lower bounds for $h(D)$ has not received nearly as much attention as finding upper bounds.





\end{document}